\documentclass[10pt, a4paper]{amsart}
\input xy   
\xyoption{all}
\numberwithin{equation}{section}

\swapnumbers

\theoremstyle{plain}

  \begingroup
        \newtheorem{theorem}[equation]{Theorem}
        \newtheorem{lemma}[equation]{Lemma}
        \newtheorem{proposition}[equation]{Proposition}
        \newtheorem{corollary}[equation]{Corollary}

        \newtheorem{remark}[equation]{Remark}
        \newtheorem{definition}[equation]{Definition}
        \newtheorem{notation}[equation]{Notation}
  \endgroup
\theoremstyle{definition}
  \begingroup
        \newtheorem{example}[equation]{Example}
        
  \endgroup

\newcommand{\mr}[1]{\buildrel {#1} \over \longrightarrow}

\newcommand{\siff}{\Leftrightarrow}

\newcommand{\rimply}{\Rightarrow}

\newcommand{\uimply}{\Uparrow}

\newcommand{\cc}{\mathcal}
\newcommand{\bb}{\mathbb}

\usepackage{colordvi}

\begin{document}

\title{Topological functors as familiarly-fibrations}

\author{Eduardo J. Dubuc and Luis Espa\~nol}

\vspace{2ex}

\begin{abstract}

In this paper we develop the theory of topological categories
over a base category, that is, a theory of topological functors. Our notion of topological
functor is similar to (but not the same) the existing notions in the
literature (see \cite{BO} 7.3), and it aims at the same examples. 
In our sense, a (pre) topological functor is a functor that creates cartesian
families. A topological functor is, in particular, a fibration, and our emphasis is put in this fact.

\end{abstract}

\maketitle

\vspace{8ex} 

{\sc introduction} \indent 
In this paper we develop the theory of topological categories over a 
base category, that is, a theory of topological functors.   Our notion of topological
functor is similar to (but not the same) the existing notions in the
literature (see \cite{BO} 7.3), and it aims at the same examples. 

Recall that a (pre) fibration
is  a functor that creates cartesian arrows. In our sense,
a pre-topological functor is a functor that creates cartesian
families, and it is topological provided that these families
compose. A topological functor is, in particular, a fibration, and
our emphasis is put in this fact.  We develop an adequate
generalization utilizing cartesian families (instead of cartesian
arrows) of the basic ideas of Grothendieck's theory of fibered categories. 

\vspace{1ex}
 
In section \ref{families} we set the basic facts of a systematic theory
of families of arrows in a category. In section \ref{se_vs_fs}
we consider \emph{u-final} and
\emph{u-surjective} families in a category $\bb{T}$ with respect to
a functor $\bb{T} \mr{u} \bb{S}$, and prove a
general theorem that characterizes the intrinsic \emph{strict
  epimorphic families in $\bb{T}$ } as the $u$-final and $u$-surjective
families.  This theorem proves to be very useful in practice, allowing
to generalize to a general setting many of the usual arguments and
constructions known
for the category of topological spaces. The assumption on $u$ in this
theorem defines the concept of \mbox{$\cc{E}$-functor,} which determines the
right generality for these constructions.
In \mbox{section \ref{cartfam}}  
we develop the basic yoga of 
cartesian arrows introduced by Grothendieck in \cite{G1}, but we do
so using families instead of single arrows.
In particular, given a fibration, we clarify the relation between
cartesian families 
and initial families, and, related to this, the relation of these
families with products in the fibers. In section \ref{Topfuns} we
define topological functors and
prove in our context all the properties corresponding to the usual
properties of topological
functors, and several new characterizations of these functors. In
particular, a characterization of topological functors in terms
of the pseudofunctor associated to the fibration.

\vspace{1ex}

For all the  concepts considered in this paper  there is a
  corresponding dual concept, and  all the corresponding dual
  statements (dual assumption and dual conclusion)  hold. We will explicitly dualize concepts and statements only when it is 
convenient or necessary.  

\pagebreak

\tableofcontents

\section{Families of arrows}  \label{families}

In this section we recall some notions and results that we shall explicitly need
in the following sections, and in this way fix notation and
terminology. 

Given a category $\mathbb{T}$ and an object $X$ in $\mathbb{T}$, we 
shall work with families 
\mbox{$(X_{\alpha} \mr{g_{\alpha}} X)_{\alpha \in \Gamma}$} of 
arrows of $\mathbb{T}$ with codomain $X$. 
Dually, we can also consider families of arrows with domain $X$.

\begin{notation}
Given a family $(X_{\alpha} \mr{g_{\alpha}} X)_{\alpha \in \Gamma}$,
    we shall simply write $X_{\alpha} \mr{g_{\alpha}} X$, 
    omitting as well a label for the index set (the context will
    always tell whether we are considering a single $\alpha$ or the
    whole family). 

The diagrammatic notation always denotes a commutative diagram, unless
otherwise explicitly indicated.
\end{notation}
We denote $\bb{T}(-, X)$ the family of all arrows $Y\to X$ of 
$\mathbb{T}$. It is important to point 
out that we allow the families to be \emph{large}, that is, not
indexed by a set.

 Recall that a \emph{crible} is a sub-family $P \subseteq \bb{T}(-,
 X)$ such that any composite \mbox{$Z\to Y\to X$} belongs to $P$ if
 $Y\to X$ belongs to $P$.

\vspace{1ex}

\begin{definition} \label{refinement}
We say that a family $Y_\lambda  \to X$ \emph{refines (is a refinement
 of)}  a
family $X_\alpha \to  X$ if there is a function between the indices 
$\lambda \mapsto \alpha_\lambda$ together with arrows 
$Y_\lambda \to X_{\alpha_\lambda}$ such that 
$$
\xymatrix@R=15pt@C=7pt
                 {
                  Y_\lambda  \ar[rr]  \ar[rd] &&  X_{\alpha_\lambda} \ar[ld]
                  \\
                  & X
                 }
$$
\end{definition} 

Any family $X_\alpha \to  X$ is a refinement (but not a sub-family) of
$\bb{T}(-, X)$ in a canonical way. We denote $P(X_\alpha \to  X)$
the crible of all arrows $Y\to X$ factorizing by some arrow $X_\alpha
\to  X$. The family $X_\alpha \to  X$ refines in a canonical way the
family $P(X_\alpha \to  X)$, and in the other direction,  $P(X_\alpha
\to  X)$ also refines  $X_\alpha \to  X$, but there is no canonical
refinement.

\begin{definition} \label{rpullback}
Given an arrow $Y \to X$, we say than a family $Y_\lambda  \to Y$ is a 
\mbox{$r$\emph{-pull-back}} of a family  $X_\alpha \to  X$, if 
$Y_\lambda  \to Y \to X$ refines $X_\alpha \to  X$. That is, if there is a function between the indices $\lambda \mapsto \alpha_\lambda$ together with arrows 
$Y_\lambda \to X_{\alpha_\lambda}$ such that 
$$
\xymatrix@R=15pt@C=20pt
        {
         Y_{\lambda} \ar[r] \ar[d] 
        & X_{\alpha_{\lambda}} \ar[d]
        \\
        Y \ar[r] 
        & X
        }
$$
\end{definition}

\begin{remark} \label{pullbackcrible}
Given an arrow $Y \to X$, among the r-pullbacks there is a largest one, namely, the pullback crible $P \subset \bb{T}(-, Y)$ defined by: 
$$
\xymatrix@R=5pt
              {
               \\
                Z \to Y  \in  P  \iff  \exists \alpha \;\; and\;\; Z \to X_\alpha \; \;such\;that
                \;\;\;\;\; 
               } 
\xymatrix@R=15pt@C=20pt
        {
         Z \ar[r] \ar[d] 
        & X_{\alpha} \ar[d]
        \\
        Y \ar[r] 
        & X
        }
$$
\enlargethispage*{100pt}
All the r-pullbacks are refinements of $P$.
\end{remark}
\pagebreak 

We consider collections $\cc{A}$ of classes of 
families of arrows with common codomain, one class $\cc{A}_X$ 
(eventually empty) for each object $X$ in $\mathbb{T}$. We say that a
 family in $\cc{A}_X$ is a $\cc{A}$-family over $X$. 

\vspace{1ex}

Properties of families of arrows determine a corresponding
collection. For instance, epimorphic families define a collection
$\cc{E}pi$ by:
$X_{\alpha} \mr{g_\alpha} X \,\in \, \cc{E}pi_X$ if given arrows 
$f, f':X\to Y$, the condition ($f\circ g_{\alpha}=f'\circ g_{\alpha}$
for all $\alpha$) implies $f=f'$.

\vspace{1ex}

We now define some operations on
these collections that yield new collections out of given ones:

\begin{definition}[operations on collections] \label{familyoperations} 
$ $ 

(1) We denote by $\cc{I}so$ the collection whose only arrows are the
isomorphisms.

(2) Given two collections $\cc{A}$, $\cc{B}$ we  define the \emph{composite} 
$ \cc{C} = \cc{A} \circ \cc{B}$ by means of the following implication:  
$$
  X_{\alpha} \rightarrow X \; \in \; \cc{A}_X \;\; and \;\;
  \forall \, \alpha \;\; X_{\alpha,\, \beta} \rightarrow X_{\alpha}
\; \in \; 
  \cc{B}_{X_\alpha}  
  \Longrightarrow \;\;
  X_{\alpha, \, \beta} \rightarrow X_{\alpha} \rightarrow X \;
  \in \; \cc{C}_X
$$
 
(3) Given $\cc{A}$ we define a new collection, denoted $\pi \cc{A}$, by: 
\vspace{1ex}

$Y_{\alpha} \rightarrow Y \in \pi \cc{A}\;\; \iff \;\;$ there is
$X_\alpha \to X  \in \cc{A}$ and $\;Y\to X$ such that:
$$
\xymatrix@R=5pt{\\ the \; squares\;\;}
\xymatrix@R=15pt@C=20pt
              {
                Y_\alpha \ar[r] \ar[d] 
                & X_\alpha \ar[d]
                \\
                Y \ar[r] 
               & X
              }
\xymatrix@R=5pt{\\ \;\; are \; pullbacks \; for \; all \; \alpha.}
$$

(4) Given $\cc{A}$ we define a new collection, denoted $s \cc{A}$, by:

\vspace{1ex}

$X_\alpha  \to X \;\in\; s \cc{A} \;\; \iff \;\;$ there is a
refinement by a family $Y_\lambda \to X \; \in \cc{A}_X$.
\end{definition}                        

Notice that $\cc{A}\,\subseteq\, \pi \cc{A}$, and $\cc{A}\,\subseteq\, s \cc{A}$. We set now some properties of collections $\cc{A}$ defined by means of these operations:

\begin{definition}[properties of collections] \label{familyproperties}
$ $   
 \vspace{1ex}
 
(I) \emph{Isomorphisms}:  $\cc{I}so \,\subseteq \, \cc{A}.$ 

\vspace{1ex}

(C) \emph{Closed under composition}: 
                            $\cc{A} \circ \cc{A} \,\subseteq\, \cc{A}.$

\vspace{1ex}

(U) \emph{Universal}: Given $X_\alpha \to X \, \in \cc{A}$ and $Y \to X$, there exists an r-pull-back $Y_\lambda \to Y \, \in \cc{A}$:
$$
\xymatrix@R=15pt@C=20pt
        {
         Y_{\lambda} \ar[r] \ar[d] 
        & X_{\alpha_{\lambda}} \ar[d]
        \\
        Y \ar[r] 
        & X
        }
$$

That is, $P \in s\cc{A}_Y$, where $P$ is the pullback crible in remark \ref{pullbackcrible}.

\vspace{1ex}

(S) \emph{Saturated}: $s \cc{A} \,\subseteq\, \cc{A}$ (hence $\cc{A}=s \cc{A}$).

\vspace{1ex}

(F) \emph{Filtered}: Given 
$X_{\alpha}\to X \, \in \, \cc{A}_X,\;  Y_{\beta}\to X \, \in \,
  \cc{A}_X$, there exists a common refinement 
$Z_{\lambda}\to X \; \in \; \cc{A}_X$:  
$$
\xymatrix@R=15pt@C=20pt
         {
             Z_{\lambda}  \ar[d]  \ar[r] \ar[rd]
            &  X_{\alpha_{\lambda}} \ar[d] 
          \\
          Y_{\beta_{\lambda}} \ar[r]  &  X 
         }
$$ 
\end{definition}

Notice that the collection $\cc{E}pi$ satisfies (I),  (C) and (S), but fails in general to satisfy (U) and (F).

\begin{proposition} \label{UCimpliesF}  
The following statements hold:

(i) If a collection satisfies the properties (U) and (C) then it satisfies (F).

(ii) If the collection $\cc{A}$ satisfies (S), then (U) is equivalent to the condition that says that 
$P \in \cc{A}_Y$, where $P$ is the pullback crible in remark \ref{pullbackcrible}.

(iii) If  the collection $\cc{A}$ satisfies (S), and the category has finite limits, then (U) is 
equivalent to the condition of stability under pulling-back, that is, $\pi\cc{A} \subseteq \cc{A}$ (hence $\cc{A}=\pi\cc{A}$).
\qed
\end{proposition}

In practice some collections are determined by the conjunction of two properties, therefore they are the intersection of two different collections. Concerning this we have: 

\begin{remark} \label{intersection}
Given two collections $\cc{A}$, $\cc{B}$, consider the collection $\cc{A} \cap \cc{B}$. Then:

i) If $\cc{A}$ and $\cc{B}$ both satisfy (I), (resp. (C)), (resp. (S)), then so does $\cc{A} \cap \cc{B}$.
 
This is not the case for conditions (U) and (F). However, it will be so for saturated collections (argue with the pullback crible):
 
ii) If  $\cc{A}$ and $\cc{B}$ both satisfy (S) and (U), (resp. (S) and (F)), then so does  $\cc{A} \cap \cc{B}$.
\end{remark} 

An important collection of families are the strict epimorphic
families. We recall now this notion from SGA4 \cite[I, 10.3, p. 180]{G2}:  

\begin{definition} \label{strictepi}
Given two families of arrows 
$f_{\alpha}:X_{\alpha}\to X$,  $g_{\alpha}:X_{\alpha}\to Y$, 
with the same indexes and domains, we say that  $g_{\alpha}$ is 
\emph{compatible with} $f_{\alpha}$ if for any pair of arrows 
$(x_{\alpha}:Z\to X_{\alpha}, \; x_{\beta}:Z\to X_{\beta})$ with the 
same domain the following condition holds: 
$f_{\alpha}\circ x_{\alpha}= f_{\beta}\circ x_{\beta}$ implies 
$g_{\alpha}\circ x_{\alpha}= g_{\beta}\circ x_{\beta}$.
 
A family $f_{\alpha}:X_{\alpha}\to X$ is \emph{strict epimorphic} 
if for any family $g_{\alpha}:X_{\alpha}\to Y$ \mbox{which is} 
compatible with $f_{\alpha}$, there exists a unique $g:X\to Y$ 
such that $g\circ f_{\alpha}=g_{\alpha}$ for all $\alpha$. 
\end{definition}

The situation is described in the following diagram, where the 
family $g_{\alpha}$ is compatible with the family $f_{\alpha}$:
$$        
\xymatrix@1
        {
         X_{\alpha}\;\; \ar @<+2pt> `u[r] `[rr]^{g_{\alpha}} [rr]
                                             \ar[r]^{f_{\alpha}} 
         & \;\;X\;\;  \ar@{-->}[r]^{\exists ! g} 
         & \;\;Y      
        }
$$

In $Set$ the category of sets (as in any topos), every epimorphic family is strict epimorphic
and the collection satisfies all five properties in definition 
\ref{familyproperties}. It is an important exactness property for a
category when the collection of strict epimorphic families satisfy all
five properties in definition \ref{familyproperties}. This is the case for the regular categories. 

It follows from the uniqueness that strict epimorphic families are epimorphic, but the converse is not true in general. Recall that the following holds:

\begin{remark}  \label{colimitconearestrictepi}
The  family of all the inclusions of a colimit cone is a strict epimorphic family.
 \qed  \end{remark}

Let $s\cc{E}_X$ be the class of all strict 
epimorphic families with codomain $X$. The collection 
$s\cc{E}$ satisfies (I) and (S), but in general it fails to satisfy 
(C), (U) and (F) in definition \ref{familyproperties}. 
Recall the following:

\begin{proposition} \label{preservation1}
If a functor has a right adjoint, then it
preserves $\cc{E}pi$-families and $s\cc{E}$-families (dually, we have
that monomorphic, and strict monomorphic 
families are preserved by functors with a left adjoint).
\qed \end{proposition}

The following is a technical fact that we shall need later, and it is
easy to prove: 
\begin{proposition} \label{compatibility} 
 If a functor has a left adjoint, then it preserves the relation of 
compatibility between families.
\qed \end{proposition}

An strict epimorphisms is the coequalizer of its kernel pair (if the
later exists). In this case some authors call them regular
epimorphisms. Strict epimorphisms rather than epimorphisms are
relevant in the important cases of factorizations.
 
\vspace{2ex}

\section{Strict epimorphic and final surjective families}
\label{se_vs_fs}  
\vspace{1ex}

In $Top$, the category of topological spaces and continuous maps, a
family is strict epimorphic  if and only if it is surjective and the
codomain has the final topology.  In this section we develop a
similar characterization of strict epimorphic families in a general
setting.  

We consider a category $\mathbb{T}$ endowed with a functor
$u:\mathbb{T}\to\mathbb{S}$ which  will satisfy some general
conditions, and we consider notions associated to the functor $u$ for 
families in $\mathbb{T}$.

We start with some definitions:
\begin{definition} \label{terminology}
An object $X$ in $\mathbb{T}$ \emph{sits over} an object $S$ in 
$\mathbb{S}$ when $u(X)=S$. We say also that $X$ is an object 
\emph{over} $S$. An arrow $f:X\to Y$ in $\mathbb{T}$ \emph{sits over}
an arrow $\varphi:S\to T$ in $\mathbb{S}$ when $u(f)=\varphi$, so  
that $X$ (resp. $Y$) sits over $S$ (resp. $T$). We say also that 
$\varphi$ \emph{lifts to} an arrow in  $\mathbb{T}$ when there exists
$f$ over $\varphi$. A family $f_\alpha:X\to Y$ in $\mathbb{T}$  
\emph{sits over} a family $\varphi_\alpha:S\to T$ in $\mathbb{S}$ 
when $u(f_\alpha)=\varphi_\alpha$, for any $\alpha$. We say also that
the family $\varphi_\alpha$ \emph{lifts to} a family in $\mathbb{T}$
when there exists $f_\alpha$ over  $\varphi_\alpha$.
\end{definition}
\begin{definition} 
We say that two families $f_{\alpha}:X_{\alpha}\to X$, 
$g_{\alpha}:X_{\alpha}\to Y$ in $\mathbb{T}$ which sit 
over the same family $\phi_{\alpha}:S_{\alpha}\to S$ in $\mathbb{S}$  
are u\emph{-isomorphic} if there exists an isomorphism $\theta:X\to 
Y$ over $id:S\to S$ such that $\theta\circ f_{\alpha}=g_{\alpha}$ for 
all $\alpha$.
\end{definition}
\begin{remark}
All collections considered in this paper are assumed to be closed
under u-isomorphisms without need to say so explicitly. 
\end{remark} 
\begin{definition} \label{uniqueuptoiso}
Let $\cc{A}$ be a collection of classes of families in 
$\mathbb{T}$. We say that \mbox{$\cc{A}$-families} are \emph{unique up to
isomorphisms} if given any two $\cc{A}$-families 
\mbox{$f_{\alpha}:X_{\alpha}\to X$,} $g_{\alpha}:X_{\alpha}\to Y$
which sit over the same family $\phi_{\alpha}:S_{\alpha}\to S$ in
$\mathbb{S}$, they are \mbox{u-isomorphic} by a unique isomorphism. 
\end{definition}

Notice that when $u$ is faithful strict epimorphic families are
unique up to \mbox{isomorphisms.} 

\begin{definition} \label{creation}
Consider collections $\cc{A}$ in $\mathbb{T}$ and 
$\cc{B}$ in $\mathbb{S}$. We say that $u$ \emph{creates} 
$\cc{A}$\emph{-families over} $\cc{B}$\emph{-families} if given any 
$\cc{B}$-family $\phi_{\alpha}:S_{\alpha}\to S$, and an object 
$X_{\alpha}$ over $S_{\alpha}$ for every $\alpha$, there exists an 
$\cc{A}$-family $f_{\alpha}:X_{\alpha}\to X$ over 
$\phi_{\alpha}:S_{\alpha}\to S$.

When the class $\cc{B}$ in $\bb{S}$ is the ``same'' class than the
class $\cc{A}$ in $\bb{T}$ (that is, if they are denoted by the same letter), we simply say that $u$ \emph{creates} $\cc{A}$\emph{-families}. 
\end{definition}

Complementing remark \ref{colimitconearestrictepi} it is
immediate to check:

\begin{remark} \label{preservescolandlim}
A functor which creates and preserves strict epimorphic families
 (finite strict epimorphic families)  create and preserve any
  colimits (resp. finite colimits) that may exists.
\qed\end{remark}

The following lemma is the key to the proof of the characterization
theorem \ref{characterization2}.  

\begin{lemma} \label{A=C}
Let $\cc{A}$, $\cc{C}$ be collections in $\mathbb{T}$, and $\cc{B}$ 
in $\mathbb{S}$ such that: 

$\cc{A}$-families are $\cc{C}$-families (i.e.,
$\cc{A}\subseteq\cc{C}$). 

$u$ sends $\cc{C}$-families into $\cc{B}$-families (i.e.,
$u(\cc{C})\subseteq\cc{B}$). 

$u$ creates $\cc{A}$-families over $\cc{B}$-families.

$\cc{C}$-families are unique up to isomorphisms.

\vspace{1ex}

Then $\cc{C}\subseteq\cc{A}$, thus  $\cc{C} = \cc{A}$
\end{lemma}
\begin{proof}
Let $X_\alpha \mr{f_\alpha} X$ be a $\cc{C}$-family over $S_\alpha
\mr{\varphi_\alpha} S$. It is isomorphic with the 
$\cc{A}$-family created over $S_\alpha \mr{\varphi_\alpha} S$. It
follows that it is a $\cc{A}$-family. 
\end{proof}  

\vspace{1ex} 

\begin{definition} Given a functor $u:\mathbb{T}\to\mathbb{S}$, we
say that a family $f_{\alpha}:X_{\alpha}\to X$ in $\mathbb{T}$ is  
u\emph{-surjective} when the family $u(f_{\alpha})$ is strict 
epimorphic in $\mathbb{S}$.
\end{definition}

We shall often omit the $u$ when we write ``$u$-surjective'' and
simply write 
``surjective''. Notice that when $u=id$, the surjective families are
the strict epimorphic families.

\vspace{1ex}

Let $\cc{S}_X$ be the class of all surjective 
families with codomain $X$. The collection \mbox{$\cc{S}$} 
satisfies conditions (I) and (S), but fails in general
to satisfy \mbox{(C), (U) and (F).} 

\begin{definition} \label{final}
Given a functor $u:\mathbb{T}\to\mathbb{S}$, let 
$f_{\alpha}:X_{\alpha}\to X$ be a family in  $\mathbb{T}$ over 
$\varphi_{\alpha}:S_{\alpha}\to S$ in $\mathbb{S}$. The family 
$f_{\alpha}$ is u\emph{-final} if  for any family 
$g_{\alpha}:X_{\alpha}\to Y$ in $\mathbb{T}$ and arrow $\phi:S\to 
T$ in $\mathbb{S}$  such that $g_\alpha$ sits over 
$\phi\circ\varphi_\alpha$, there exits a unique $g:X\to Y$ over 
$\phi$ such that $g\circ f_\alpha=g_\alpha$. 
\end{definition}
We shall often omit the $U$ when we write ``$u
$-final'' and simply
write ``final''.  Notice that when $u=id$ all families are final.
The situation for final families is described in the
following double diagram, where the top diagram sits over the bottom
diagram. 
\begin{equation}   \label{doublediagram}     
\xymatrix@1@R=15pt
        {
         X_{\alpha}\;\; \ar @<+2pt> `u[r] `[rr]^{g_{\alpha}} [rr]
                                             \ar[r]^{f_{\alpha}} 
         & \;\;X\;\;  \ar@{-->}[r]^{\exists ! g} 
         & \;\;Y  
         \\        
         S_{\alpha}\;\; \ar[r]^{\varphi_{\alpha}} 
         & \;\;S\;\;  \ar[r]^{\phi} 
         & \;\;T      
        }
\end{equation}
Final families are unique up to isomorphisms in the sense of
definition \ref{uniqueuptoiso}

Let $\cc{F}_X$ be the class of all final 
families with codomain $X$. The collection $\cc{F}$ 
satisfies conditions (I), (C) and (S), but fail in 
general to satisfy (U) and (F). 

By $\cc{F}\cc{S}=\cc{F} \cap \cc{S}$ we shall denote the
collection of all final and surjective families. Notice that
$\cc{FS}$-families are unique up to isomorphisms. Moreover, the collection $\cc{FS}$ 
satisfies conditions (I) and (S) by remark \ref{intersection}. 

\vspace{1ex}

Corresponding to proposition  \ref{preservation1} it is easy to prove
the following:

\begin{proposition} \label{preservation2}
Given a commutative triangle of functors, if $F$ 
has a right adjoin $R$ such that $u\circ R=u'$, then $F$ preserves 
final families.
$$
\xymatrix@R=15pt@C=7pt
         {
          \mathbb{T} \ar[rr]^{F}  \ar[rd]_{u}  
           &&  \mathbb{T}' \ar[ld]^{u'}  
           \\
           & \cc{S} 
         }
$$ 
\qed \end{proposition}
 We remark that if the right 
adjoint $R$ fails to make the triangle commutative, then $F$ will not 
preserve general final families. However, very often in practice it
will preserve those final families which are also surjective. 

\vspace{1ex}

We study now the implication  
$:\;\;Strict \; epimorphic \;\; \Rightarrow \;\; Final\; surjective$ 

\begin{proposition} \label{strictepiaresurjective} 
If $u$ preserves strict epimorphic families then by definition strict
epimorphic families are surjective, that is,  
$s\cc{E}\subseteq\cc{S}$. By Proposition \ref{preservation1}, this is
the case when  
$u$ has a right adjoint.
\qed \end{proposition}

\begin{proposition} \label{strictepiarefinal}
If $u$ is faithful, strict epimorphic families are final, that 
is, $s\cc{E}\subseteq\cc{F}$.
\end{proposition}  
\begin{proof} 
Consider the double diagram \ref{doublediagram}, and suppose that
 $f_\alpha$ is strict epimorphic. Since the family
 $\phi \circ \varphi_\alpha$ is compatible with $ \varphi_\alpha$, it follows from the 
faithfulness of $u$ that $g_\alpha$ is compatible with
$f_\alpha$. The existence of $g$ follows. 
\end{proof}

\begin{corollary} \label{strictepiarefinalsurjective}
If $u$ is faithful and preserves strict epimorphic families then each 
strict epimorphic family is final and surjective, that is, 
$s\cc{E}\subseteq\cc{FS}$.
\end{corollary}  
\begin{proof} 
Apply propositions \ref{strictepiaresurjective} and
\ref{strictepiarefinal}.  
\end{proof}

From proposition \ref{colimitconearestrictepi} it follows:
\begin{corollary}  \label{colimitconearefinal}
If $u$ is faithful, the  family of all the inclusions of a colimit cone is a final family. 
If in addition $u$ preserves strict epimorphic families, it is final surjective. \qed
\end{corollary}

\vspace{1ex}

We pass now to the implication 
$:\;\; Final \; surjective  \;\; \Rightarrow \;\; Strict \;
epimorphic$

\begin{proposition} \label{finalsurjectivearestrictepi}
If $u$ has a left adjoint, then each final surjective family is
strict epimorphic, that is $\cc{FS}\subseteq s\cc{E}$. 
\end{proposition}
\begin{proof}
Let $X_\alpha \mr{f_\alpha} X$ be a final surjective family. Given
any compatible family  $X_\alpha \mr{g_\alpha} Y$, by proposition 
\ref{compatibility}  $uX_\alpha \mr{ug_\alpha} uY$ is compatible with
$uX_\alpha \mr{uf_\alpha} uX$, which by assumption is strict
epimorphic in $\cc{S}$. Thus we have $uX \mr{\varphi} uY$. This
finishes the proof since we are assuming also that  $X_\alpha
\mr{f_\alpha} X$ is a final family. 
\end{proof}
From \ref{preservation1}, \ref{strictepiarefinalsurjective} and
\ref{finalsurjectivearestrictepi} it follows: 
 
\begin{theorem} \label{characterization1}
If the functor $u$ is faithful, and has left and right adjoints, then
a family is strict epimorphic if and only if is final  
and surjective, that is, \mbox{$s\cc{E}=\cc{FS}$.}
\qed\end{theorem}

Notice that the hypothesis in this theorem are self dual, so that a
dual theorem also holds:

\begin{theorem} \label{characterization1dual}
If the functor $u$ is faithful, and has left and right adjoints, then
a family is strict monomorphic if and only if is initial  
and injective.
\qed\end{theorem}

We have the same characterization of strict epimorphic families under
a different assumption which has also many other important
consequences. This assumption defines the right generality for the
validity of many constructions and results, and merits to be treated by itself.

\begin{definition} \label{efunctor}
A functor $\bb{T} \mr{u} \bb{S}$ is a $\cc{E}$-functor (resp.
$\cc{M}$-functor) if it is \mbox{faithful} and creates and preserves strict
epimorphic families (resp. strict monomorphic families).
\end{definition}

If we consider only finite families in definition \ref{efunctor}, we
have the notions of  \mbox{$f\cc{E}$-functor} and
$f\cc{M}$-functor. 
  
\begin{theorem}\label{characterization2}
Given a $\cc{E}$-functor  $\bb{T} \mr{u} \bb{S}$, a family in $\bb{T}$
is strict epimorphic if and only 
if is final  and surjective, that is, $s\cc{E}=\cc{FS}$.
\end{theorem}  
\begin{proof} 
By Corollary \ref{strictepiarefinalsurjective},
$s\cc{E}\subseteq\cc{FS}$, and we know that  
$\cc{FS}$-families are unique up to isomorphisms. By Lemma \ref{A=C},
it remains to see that $u$ creates $s\cc{E}$-families over
$u(\cc{FS})$-families.  
But in $\mathbb{S}$ by definition we have 
$u(\cc{FS})\subseteq u(\cc{S})\subseteq s\cc{E}$, and by assumption
$u$  creates strict epimorphic families.
\end{proof}

The reader can easily check that all the arguing above still holds if
we consider only finite families. We state now explicitly the dual
statement of theorem \ref{characterization2} in the finite case: 
\begin{theorem} \label{finitecharacterization2} 
Given a $f\cc{M}$-functor  $\bb{T} \mr{u} \bb{S}$, a finite family in
$\bb{T}$ is strict monomorphic if and only 
if is initial and injective.
\qed\end{theorem}

\vspace{2ex}

\section{Cartesian families}  \vspace{1ex}  \label{cartfam}

In this section we develop the basic yoga of 
cartesian arrows introduced by Grothendieck in \cite{G1}, but we do
so using families instead of single arrows. The consequences of
replacing arrows by families are strong and conform a very different
theory which furnishes the framework for the theory of topological functors.

Given a functor $u:\mathbb{T}\to\mathbb{S}$ and an object $S \in
\bb{S}$, we denote $\bb{T}_S$ the \emph{fiber of $u$ over $S$}, that
is, the subcategory of $\bb{T}$ formed by all arrows sitting over the
identity $S \mr{id} S$ of the object $S$. Notice that each fiber $\bb{T}_S$ is a poset if $u$ is faithful.

We explicitate now the definition of initial 
family by dualizing definition \ref{final}.

\begin{definition} \label{initial}
 Given a functor $u:\mathbb{T}\to\mathbb{S}$ and a family
$f_{\alpha}:X\to X_{\alpha}$  in  $\mathbb{T}$ over 
$\varphi_{\alpha}:S\to S_{\alpha}$ in $\mathbb{S}$, we say that the family 
$f_{\alpha}$ is \emph{$u$-initial} if  for any family 
$g_{\alpha}:Y\to X_{\alpha}$ in $\mathbb{T}$ and arrow $\varphi:R\to 
S$ in $\mathbb{S}$  such that $g_\alpha$ sits over 
$\varphi_\alpha\circ\phi$, there exits a unique $g:X\to Y$ over 
$\varphi$ such that $f_\alpha\circ g=g_\alpha$. We shall often omit the $u$ when we 
mean ``$u$-initial'' and simply
say ``initial''.
\end{definition}

The situation for initial families is described in the following 
double diagram, where the top diagram sits over the bottom diagram.
\begin{equation}    \label{doublediagram2}    
\xymatrix@1@R=15pt
        {
         Y\;\; \ar @<+2pt> `u[r] `[rr]^{g_{\alpha}} [rr]
                                             \ar@{-->}[r]^{g} 
         & \;\;X\;\;  \ar[r]^{f_{\alpha}} 
         & \;\;X_{\alpha}  
         \\  
         R\;\; \ar[r]^{\varphi} 
         & \;\;S\;\;  \ar[r]^{\varphi_{\alpha}} 
         & \;\;S_{\alpha}      
        }
\end{equation}
If we consider, for any object $X$ in $\mathbb{T}$, the class 
$\cc{I}_X$ of all initial families with domain $X$, then the
collection $\cc{I}$ satisfies (I), (C) and (S$^{op}$), but fail in
general to satisfy (U$^{op}$) and (F$^{op}$).


Notice that an initial family $f_{\alpha}$ acts as a monomorphic family on arrows $g$, $g'$ over  the same $\varphi$.

\begin{definition} \label{cartesian}
 Given a functor $u:\mathbb{T}\to\mathbb{S}$ and a family 
$f_{\alpha}:X\to X_{\alpha}$ in  $\mathbb{T}$ over 
$\varphi_{\alpha}:S\to S_{\alpha}$ in $\mathbb{S}$, we say that the family 
$f_{\alpha}$ is $u$-cartesian if  for any family 
$g_{\alpha}:Y\to X_{\alpha}$ in $\mathbb{T}$ over 
$\varphi_{\alpha}:S\to S_{\alpha}$, there exits a unique $g:Y\to 
X$ in $\bb{T}_S$ such that $f_\alpha\circ g=g_\alpha$. We shall often omit the $u$ 
when we mean ``$u$-cartesian'' and simply say ``cartesian''.
\end{definition}

The situation for cartesian families is described in the following 
double diagram, where the top diagram sits over the bottom diagram.
\begin{equation}        
\xymatrix@1@R=15pt
        {
         Y \ar@/^5pt/[rd]^{g_{\alpha}} \ar@{-->}[d]_{g} 
         \\
         X   \ar[r]^{f_{\alpha}} & \;\;X_{\alpha} 
         \\ 
         S\;\;  \ar[r]^{\varphi_{\alpha}} & \;\;S_{\alpha}      
        }
\end{equation}  

\begin{notation}
In a double diagram, the top part always sits over the bottom part,
and all vertical arrows always sit over the respective identity arrow
(that is, vertical arrows are always in the fibers).
\end{notation}

Cartesian families are unique up to isomorphisms in the sense of
definition  \ref{uniqueuptoiso}. Notice then that the isomorphism
lives in the fiber $\bb{T}_S$.

\vspace{1ex}

When the family $f_{\alpha}$ is formed by a unique arrow $X \mr{f} X'$,  
we recover Grothendiek's SGA-notion of cartesian arrow. 

We say that  functor $u$ \emph{creates}
cartesian (resp. initial) arrows if for any arrow $T\mr{\varphi} S$ in
$\mathbb{S}$ and object $X \in \bb{T}$ over $S$, there exists $Y \in
\mathbb{T}$ over $T$, and a cartesian (resp. initial) arrow $Y\mr{f}
X$ over $\varphi$.


\vspace{1ex}

Recall \cite[t.2, VI, 6.1, p. 102]{G2} the following definition:
\begin{definition}
A functor $u: \bb{T} \to \bb{S} $ is a \emph{prefibration} (equivalently,  $\bb{T}$ is 
\emph{prefibered}) if $u$ creates cartesian arrows, and it is a \emph{fibration} 
(equivalently, $\bb{T}$ is \emph{fibered}) if it is a prefibration and cartesian arrows 
compose.
\end{definition}

Recall that initial arrows compose and are cartesian, and that when $\bb{T}$ is prefibered, if cartesian arrows compose then they are initial. Thus:

\begin{proposition} \label{fibrations}
A functor $u: \bb{T} \to \bb{S} $ is a fibration  if and only if it creates initial arrows.
\qed  \end{proposition}

\vspace{1ex}

We study now this situation for families. Clearly initial families are cartesian, but the converse is not true.
Cartesian families do not compose, that is, they do not satisfy 
property (C).

We establish first a technical lemma relating products in the fibers with cartesian arrows.
\begin{lemma} \label{technicallemma}
Given a functor $u:\mathbb{T}\to\mathbb{S}$, consider the situation described in the 
following double diagram:
$$        
\xymatrix@1
        {
         X\; \ar@/^5pt/[rd]^{p_{\alpha}} \ar[d]_{\pi_\alpha} 
         \\
         Y_\alpha\;  \ar[r]^{f_\alpha} & \;\;X_\alpha
         \\ 
         S\;  \ar[r]^{\varphi_{\alpha}} & \;\;S_{\alpha}      
        }
$$

(i) Assume each $f_\alpha$ cartesian, then:  
\begin{center}
$\pi_\alpha$ is a product in $\mathbb{T}_S$  $\;\; \siff \;\; $ $p_\alpha$ is a cartesian 
family.
\end{center} 

(ii) Assume each $f_\alpha$ initial, then: 
\begin{center}
$\pi_\alpha$ is an initial family  $\;\; \siff \;\; $ $p_\alpha$ is an initial family. 
\end{center}
\end{lemma} 
\begin{proof} 
(i): We shall argue over the situation described in the following double diagram:
$$
\xymatrix@C=50pt
                {
                Y  \ar@/_20pt/[dd]^{\!h_\alpha}  \ar[d]^g  \ar@/^10pt/[rdd]^{g_\alpha} 
                \\
                X \ar[d]^{\pi_\alpha}  \ar[dr]^{\!p_\alpha}
                \\
                Y_\alpha \ar[r]^{f_\alpha}  &  X_\alpha
                \\
                S \ar[r]^{\varphi_\alpha}  &  S_\alpha
                }
$$

($\Leftarrow$): Given a family $h_\alpha$, consider the composites 
$f_\alpha \circ h_\alpha$. Then there is a unique $g$ such that 
$p_\alpha \circ g =  f_\alpha \circ h_\alpha$. But $p_\alpha = f_\alpha \circ \pi_\alpha$, so that  $ f_\alpha \circ \pi_\alpha \circ g =  f_\alpha \circ h_\alpha$.
Since each $f_\alpha$ is cartesian, it follows that 
$\pi_\alpha \circ g = h_\alpha$. Given $g'$ such that 
$\pi_\alpha \circ g' = h_\alpha$, we have
$f_\alpha \circ \pi_\alpha \circ g' =  f_\alpha \circ h_\alpha$, so that  
$p_\alpha \circ g'=  f_\alpha \circ h_\alpha$. This implies $g' = g$.

($\Rightarrow$): Given a family $g_\alpha$, since each $f_\alpha$ is cartesian, there exists a unique $h_\alpha$ such that $f_\alpha \circ h_\alpha = g_\alpha$. It follows that there exists in turn a unique $g$ such that  $\pi_\alpha \circ g_\alpha = h_\alpha$. We have $f_\alpha \circ \pi_\alpha \circ g =  f_\alpha \circ h_\alpha$, thus
$p_\alpha \circ g =  g_\alpha$. Given $g'$ such that  $p_\alpha \circ g' =  g_\alpha$, we have $f_\alpha \circ \pi_\alpha \circ g' =  f_\alpha \circ h_\alpha$. Since each $f_\alpha$ is cartesian, it follows $\pi_\alpha \circ g' =  h_\alpha$. This implies $g' = g$.

\vspace{1ex}

(ii): It follows the same lines than the preceding proof, and it is left to the reader.
\end{proof}

Considering the case in which  $f_\alpha = id_{X_\alpha}$ for all indices $\alpha$, we have from part $i)$:
\begin{corollary} \label{cartesianiffproduct}
A family $\pi_\alpha:X\to X_\alpha$ in $\mathbb{T}_S$ is a product 
in $\mathbb{T}_S$ if and only if it is a cartesian family.  
\end{corollary}

 Concerning composition of cartesian families, it suffices to consider
 composition of arrows with families as follows:
\begin{definition} Given a functor $u:\mathbb{T}\to\mathbb{S}$ we say 
that \emph{cartesian arrows compose with cartesian families} if for 
any cartesian arrow $f:Y\to X$ and cartesian family 
$f_{\alpha}:X\to X_{\alpha}$, the family $f_{\alpha}\circ 
f:Y\to X_{\alpha}$ is cartesian.
\end{definition} 

\begin{proposition} \label{cartesian=initial}
Given a prefibration $u:\mathbb{T}\to\mathbb{S}$,  the following statements are 
equivalent:

(i) Cartesian arrows compose with cartesian families.

(ii) Cartesian families are initial families.

(iii) Cartesian families compose.   
\end{proposition}  
\begin{proof} 
Since initial families compose, it remains to prove the implication $(i) \Rightarrow (ii)$. Let  
$X \mr{f_\alpha} X_\alpha$ be a cartesian family over
$S_\alpha$. Consider $g_\alpha$ and $\varphi$ as in the  following double diagram:
$$
\xymatrix
                {
                 Z  \ar@{-->}[d]^h   \ar@{-->}[dr]^{\!\!\ell} \ar@/^5pt/ [drr]^{g_\alpha}
                 \\
                 Y \ar@{-->}[r]^g & X \ar[r]^{f_\alpha}  &  X_\alpha
                 \\
                 H \ar[r]^\varphi & S  \ar[r]^{\varphi_\alpha} & S_\alpha
                 }
$$
Take $g$ cartesian over $\varphi$. Since $f_\alpha \circ g$ is cartesian, we have a 
unique $h$ in $\bb{T}_H$ such that $f_\alpha \circ g \circ h = g_{\alpha}$. Let $\ell = g
\circ h$. We have that
$\ell$ sits over $\varphi$ since $h$ sits over $id_H$, and
$f_\alpha \circ \ell = g_\alpha$ since $f_\alpha \circ \ell = f_\alpha \circ g \circ h = g_
\alpha$. We leave the reader to check the uniqueness of $\ell$.
\end{proof}

We conclude then that in a prefibration \emph{cartesian families compose if and only if they are initial families.}

\begin{definition} \label{productsstable}
Given any functor $u:\mathbb{T}\to\mathbb{S}$, we say 
that \emph{products (in the fibers) are stable} if given 
$\varphi:S\to T$ in $\mathbb{S}$ and a double diagram as follows:
$$        
\xymatrix@1
        {
         Y\; \ar[r]^{g} \ar[d]^{\rho_\alpha}  &  X \ar[d]^{\pi_\alpha}
         \\
         Y_\alpha\;  \ar[r]^{g_\alpha} & \;\;X_\alpha
         \\ 
         S\;  \ar[r]^{\varphi} & \;\;T      
        }
$$
with $g$ and each arrow $g_\alpha$ cartesian, then: 
$$\; if \;\pi_\alpha \; is\; a \;product \; (in \; \mathbb{T}_T), \;so \;it \;is \;\rho_\alpha\; (in \;\mathbb{T}_S)$$
\end{definition}

\begin{proposition} \label{stableproduct=initial}
Given a prefibration $u:\mathbb{T}\to\mathbb{S}$, the following statements are equivalent:

(i) Cartesian arrows compose with cartesian families.

(ii) Cartesian arrows compose (that is, $u$ is a fibration) and products (in the fibers) are stable.   
\end{proposition} 
\begin{proof}
$(ii) \Rightarrow (i)$. Given a cartesian arrow $Y \mr{g} X$ and a cartesian family 
\mbox{$X \mr{g} Z_\alpha$} in $\bb{T}$, we shall argue over the double diagram which follows: 
$$
\xymatrix
              {
                   Y \ar[r]^g  \ar@{-->}[d] ^{\rho_\alpha} 
               &  X \ar@/^2pt/ [rd]^{p_\alpha} \ar@{-->}[d] ^{\pi_\alpha} 
               \\
                   Y_\alpha \ar@{-->}[r]^{g_\alpha}  
               &  X_\alpha  \ar@{-->}[r] ^{f_\alpha} 
               &  Z_\alpha
               \\
                   S  \ar[r]^{\varphi}  
               &  T  \ar[r] ^{\phi} 
               &  Z
              }
$$
Take $f_\alpha$ and $g_\alpha$ cartesian for each $\alpha$. It follows that there exists the vertical arrows $\pi_\alpha$ and $\rho_\alpha$. By lemma \ref{technicallemma} we have that $\pi_\alpha$ is a product. Since products are stable, so it is $\rho_\alpha$. But $f_\alpha \circ g_\alpha$ is cartesian, thus again by lemma \ref{technicallemma} we have that $p_\alpha \circ g$ is a cartesian family.

\hspace{2ex} $(i) \Rightarrow (ii)$. Consider an square as in definition \ref{productsstable}, with $\pi_\alpha$ a product. By corollary \ref{cartesianiffproduct} we have that $\pi_\alpha$ is a cartesian family, thus so it is $\pi_\alpha \circ g$. Then by lemma \ref{technicallemma} we have that $\rho_\alpha$ is  product. 
\end{proof}

\vspace{1ex}

In this paper we shall deal with prefibrations whose fibers are posets. Concerning this we have:
\begin{proposition} \label{faithful=poset}
Given a prefibration $u:\mathbb{T}\to\mathbb{S}$,  $u$ is faithful if and only if $
\mathbb{T}_S$ is a 
poset for any object $S$ in $\mathbb{S}$.   
\end{proposition}  
\begin{proof} 
Clearly if $u$ is faithful $\mathbb{T}_S$ is a poset. Conversely, given $u,v:Y\to X$ over 
$\varphi:T\to S$, take a cartesian arrow $f:Z\to X$ over $\varphi$; 
then there exist $u',v':Y\to Z$ in $\mathbb{T}_S$ such that $f\circ 
u'=u, f\circ v'=v$. But $u'=v'$, thus $u=v$.   
\end{proof}

\section{Topological functors} \vspace{1ex}  \label{Topfuns}

A \emph{(pre) topological} functor behaves as a (pre) fibration, but
with respect to families instead of single arrows. In this section we
consider four notions on a functor related as follows:

$$
\xymatrix@R=4pt@C=4pt
         {
          Pretopological & \rimply & Prefibration
          \\
          \uimply & & \uimply
          \\
          Topological & \rimply & Fibration
          }
$$ 
 
 \vspace{1ex}
 
We say that a functor $u$ \emph{creates} cartesian (resp. initial) 
families if for any  family $S \mr{\varphi_\alpha} S_\alpha$ and objects $X_\alpha \in \bb{T}$ over $S_\alpha$, there exists 
$Y \in \mathbb{T}$ over $S$ and a cartesian (resp. initial) family $Y
\mr{f_\alpha} X_\alpha$ over $S \mr{\varphi_\alpha} S_\alpha$. Notice
that in the terms of the dual case of definition \ref{creation}, this means that $u$
creates cartesian (resp. initial) families over the class
$\bb{S}(S,\,-)$ of all families in $\bb{S}$.
 
\begin{definition} \label{pretopological}
A functor $u:\mathbb{T}\to\mathbb{S}$ is called pretopological if it creates cartesian families. 
\end{definition}

From corollary \ref{cartesianiffproduct}, it follows that the fibers of a pretopological functor have all products. This property characterize such functors as prefibrations.

\begin{proposition} \label{preproduct}
A functor $u:\mathbb{T}\to\mathbb{S}$ is pretopological if and only if it is a prefibration such that the fibers  $\mathbb{T}_S$ have all products. 
\end{proposition} 
\begin{proof}
Given $X_\alpha$ over $S_\alpha$ and $S \mr{\varphi_\alpha} S_\alpha$, we shall argue over the following diagram:
$$
\xymatrix@1
        {
         X \ar@{-->}@/^5pt/[rd]^{p_{\alpha}} \ar@{-->}[d]_{\pi_\alpha} 
         \\
         Y_\alpha   \ar@{-->}[r]^{f_{\alpha}} & \;\;X_{\alpha} 
         \\ 
         S\;\;  \ar[r]^{\varphi_{\alpha}} & \;\;S_{\alpha}      
        }
$$
Take $f_\alpha$ cartesian for each $\alpha$, let $\pi_\alpha$ be a product in $\bb{T}_S$, and set $p_\alpha = f_\alpha \circ \pi_\alpha$. By lemma \ref{technicallemma}(i), it follows that $p_\alpha$ is a cartesian family.
\end{proof}

Notice that the proof of Proposition \ref{preproduct} also proves that each fiber
has small products if and only if the functor creates small cartesian
families, a fact that can be of independent interest although we will not have a 
use for it in this paper. 

As it is well known (Freyd \cite{F}) and easy to show, a category with all products collapses into a complete lattice (eventually large), thus we have

\begin{proposition} \label{preproductbis}
A functor $u:\mathbb{T}\to\mathbb{S}$ is pretopological if and only if it is a prefibration such that the fibers  $\mathbb{T}_S$ are complete lattices. 
\qed \end{proposition} 

It follows then from proposition \ref{faithful=poset} that \emph{pretopological functors are necessarially faithful}. We give now a direct proof of this fact.

\begin{proposition} [compare \cite{BO}, 7.3.4]\label{topologicalisfaithful}
Any pretopological functor is faithful.
\end{proposition}  
\begin{proof}
Let $u:\mathbb{T}\to\mathbb{S}$ be a pretopological functor. Given $f,g:X\longrightarrow Y$ over $S \mr{\varphi} T$, we shall show that $f = g$. Consider all arrows $Y \mr{\alpha}  (-)$ with domain $Y$ in $\bb{T}_S$, and for each such $\alpha$ an arrow $S \mr{\varphi_\alpha} T$, $\varphi_\alpha = \varphi$. This determines the following double diagram:
$$        
\xymatrix@1
        {
         Y \ar@/^5pt/[rd]^{h_{\alpha}} \ar@{-->}[d]_{h} 
         \\
         A   \ar[r]^{f_{\alpha}} & \;\;X 
         \\ 
         S\;\;  \ar[r]^{\varphi_{\alpha}} & \;\;T     
        }
$$
where $f_\alpha$ is a cartesian family over $\varphi_\alpha$, and $h_\alpha$ is defined as follows: 
$$
h_\alpha= 
\left \{
\begin{array}{ll}
f  &  \mbox{$if \;\; codomain(\alpha) = A \;and\; f_\alpha\circ\alpha=g$} \\
g  &  \mbox{$otherwise$}
\end{array}
\right.    
$$
It follows there exists (a unique) $h$ such that $f_\alpha \circ h = h_\alpha$ for any $\alpha$. In particular, $f_h \circ h = h_h$. We have then:
$$
h_h =
\left \{
\begin{array}{ll}
f  &  \mbox{$if \;\; f_h \circ h  = g \;\; i.e. \;\; h_h = g$} \\
g  &  \mbox{$if \;\; f_h \circ h  = f \;\; i.e. \;\; h_h = f$}
\end{array}
\right.    
$$
Since $h_h = f$ or $h_h = g$, it follows $f = g$, both being equal to $h_h$.
\end{proof}

\begin{definition} \label{topological}
A functor $u:\mathbb{T}\to\mathbb{S}$ is called 
\emph{topological} if it creates cartesian families and cartesian families compose. 
\end{definition} 

As an immediate consequence of proposition
\ref{cartesian=initial} we have the analogous result to proposition
\ref{fibrations}.
\begin{proposition} \label{topin}
A functor $u: \bb{T} \to \bb{S} $ is topological if and only if it creates initial families.
\qed  \end{proposition}

 By definition an injective family in $\bb{T}$ is a family over an strict monomorphic family in $\bb{S}$. Thus, creation of initial injective families means creation of initial families over strict monomorphic families.
A topological functor, in particular, creates initial arrows and
initial injective families. We see next that these two particular
instances of creation imply the creation of arbitrary initial families. 
\begin{proposition} \label{chartop1}
A functor $u:\mathbb{T}\to\mathbb{S}$ is topological if and only if it is a fibration (that is, it creates initial arrows) and it creates initial injective families. 
\end{proposition} 
\begin{proof}
Given $X_\alpha$ over $S_\alpha$ and $S \mr{\varphi_\alpha} S_\alpha$, we shall argue over the following diagram:
$$
\xymatrix@1
        {
         X \ar@{-->}@/^5pt/[rd]^{p_{\alpha}} \ar@{-->}[d]_{\pi_\alpha} 
         \\
         Y_\alpha   \ar@{-->}[r]^{f_{\alpha}} & \;\;X_{\alpha} 
         \\ 
         S\;\;  \ar[r]^{\varphi_{\alpha}} & \;\;S_{\alpha}      
        }
$$
Take an initial arrow $f_\alpha$ for each $\alpha$.  Consider the family of
identity arrows $id_S$, one for each $\alpha$ (obviously a strict
monomorphic family in $\bb{S}$), let $\pi_\alpha$ be an initial
family, and set $p_\alpha = f_\alpha \circ \pi_\alpha$. By lemma \ref{technicallemma}(ii), it follows that $p_\alpha$ is an initial family.
.
\end{proof}

Clearly, a topological functor is a fibration, and as an immediate corollary of propositions \ref{stableproduct=initial} and \ref{preproductbis} we have the following characterization of topological functors as fibrations:

\begin{proposition} \label{preproductbisbis}
A functor $u:\mathbb{T}\to\mathbb{S}$ is topological if and only if it is a fibration such that the fibers  $\mathbb{T}_S$ are complete lattices and infima are stable. 
\qed \end{proposition} 

We state now a characterization of topological functors $u:\mathbb{T}\to\mathbb{S}$ in terms of the pseudofunctor $P: \bb{S}^{op} \to \bb{C}at$ associated to the fibration. Recall that the data which defines a fibration $u$ is equivalent to the data that defines its associated pseudofunctor  $P$ (Grothendieck's construction, \cite{G1}). Recall that standard notation for the action of $P$ on arrows is an upper star. In this way, given $S \mr{\varphi} T$ in $\bb{S}$,  
$P(S \mr{\varphi} T) = \bb{T}_T \mr{\varphi^*} \bb{T}_S$.  Recall also that given $X$ in $\bb{T}_T$, $\varphi^{*}X$ in $\bb{T}_S$ is characterized by a cartesian arrow $\varphi^{*}X  \to X$ in $\bb{T}$ over $S \mr{\varphi} T$.

An immediate consequence of proposition \ref{preproductbisbis} is the following:
\begin{theorem}
Given a pseudofunctor $P: \bb{S}^{op} \to \bb{C}at$, the associated fibration 
$u:\mathbb{T}\to\mathbb{S}$ is a topological functor if and only if: \\
(i) For all objects $S$ in $\bb{S}$, $P(S)$ is a complete lattice. \\
(ii) For all arrows $\varphi$ in $\bb{S}$, $\varphi^*$ preserves infima, or, equivalently, it has a left adjoint $\varphi_! \dashv  \varphi^*$.
\qed \end{theorem}

We introduce now some background and notations.

\begin{notation} 
Given a functor  $u:\mathbb{T}\to\mathbb{S}$ and an object $S \in \bb{S}$,

1)  We denote by $S_\top \in \bb{T}$ the initial family over the empty family with domain $S$. If it exists, $S_\top$ is a terminal object of the fiber $\bb{T}_S$, and the assignment $S \mapsto S_\top$ furnishes a right adjoint $(-)_\top \vdash u$. 

2) We denote by $S_\bot \in \bb{T}$ the final family over the empty family with codomain $S$. If it exists, $S_\bot$ is an initial object of the fiber $\bb{T}_S$, and the assignment $S \mapsto S_\bot$ furnishes a left adjoint $u \vdash (-)_\bot$. 
\end{notation}

Notice that a right adjoint  $r:\mathbb{S}\to\mathbb{T}$ of $u$ not necessarially furnishes an initial family over the empty family. That is, $r(S)$ will not be $S_\top$. We have $r(S)=S_\top$ when the counit of the adjunction is the identity arrow, that is, essentially, when $r$ is full and faithful. Same considerations apply to a left adjoint for $u$ and $S_\bot$. 

If $1$ is a terminal object of $\mathbb{S}$, then $1_\top$ is a terminal object of $\mathbb{T}$, and  if $0$ is an initial object of $\mathbb{S}$, then $0_\bot$ is an initial object of $\mathbb{T}$.

We warn the reader that, contrary than in the usual examples, the
objects  $1_\top$  and $1_\bot$ may be different. That is, the fiber
over $1$, even in the case in which the category $\bb{S}$ is the
category of sets, will not in general be the singleton category (there
may be many different "structures" on the singleton set). 
\begin{example}
It is clear that the forgetful functor $u:Top\to Set$ is topological,
and $1$ has a unique topology (discrete $=$ indiscrete). If we
consider each set endowed with a filter of subsets instead of a
topology, and the same definition of morphisms, then we obtain a
category $Filt$ and a forgetful functor $u:Filt\to Set$ which is
topological (initial families like in $Top$) but now $1$ has two
different structures (discrete $\neq$ indiscrete). This type of
structure apears when considering convenient categories for proper homotopy theory (categories of exterior spaces, see \cite{GGH}).

\qed \end{example}


Since a topological functor creates, in particular, the empty initial
family, there exists $(-)_\top$. But also there exists $(-)_\bot$:

\begin{proposition} \label{initialisfinal}
Given a faithful functor $u:\mathbb{T}\to\mathbb{S}$, let  $X \mr{\pi_\alpha} Z$ be an initial family over the family of all arrows with domain $S$, $S \mr{\alpha} u(Z)$, $Z \in \bb{T}$. Then $X$ is the the final family over the empty family with codomain $S$,  
$X = S_\bot$ (and $\pi_\alpha$ is the arrow that corresponds to $\alpha$ by adjointness).
\end{proposition}
\begin{proof}
We have to check that $X$ is the value on $S$ of a left adjoint to $u$. This is immediate, given  $S \mr{\alpha} u(Z)$, just take $X \mr{\pi_\alpha} Z$. Uniqueness follows from the faithfulness of $u$.
\end{proof}

\begin{corollary} [compare \cite{BO} 7.3.7] \label{topologicalleftright}
Any topological functor $u:\mathbb{T}\to\mathbb{S}$ has full and
faithful left and right adjoints
$(-)_\top \vdash u \vdash (-)_\bot$, $u(-)_\top = u(-)_\bot = id$.
\qed \end{corollary}

 From proposition \ref{topologicalisfaithful} and theorems
 \ref{characterization1},  \ref{characterization1dual}, we have:
\begin{corollary} \label{fs=seintop}
Given any topological functor $u:\mathbb{T}\to\mathbb{S}$,   
a family in $\bb{T}$ is strict epimorphic if and only if it is final  
and surjective, and it is strict monomorphic if and only if it is
initial and injective.
\qed\end{corollary}

The lack of size limitations on the families not only has as a consequence the faithfulness of pretopological functors, but also implies that the notion of pretopological and topological functor is selfdual. Proposition \ref{initialisfinal} generalizes, that is, if a functor creates initial families, then it creates final families (obviously an implication equivalent to its reciproque)

\begin{proposition}[compare \cite{BO} 7.3.6] \label{topselfdual}
A functor $u:\mathbb{T}\to\mathbb{S}$ is topological if and only if considered as a functor $u:\mathbb{T}^{op} \to\mathbb{S}^{op}$ is topological.
\end{proposition}  
\begin{proof} 
It is enough to prove that if $u$ creates initial families then it creates final families.  Given a family $S_\alpha\mr{\varphi_\alpha} S$ in  $\mathbb{S}$ and a family $X_\alpha$ of objects of $\mathbb{T}$ over $S_\alpha$, we must create a final family 
$X_\alpha\mr{f_\alpha} X$ over $\varphi_\alpha$. We shall argue over the following double diagram:
$$
\xymatrix@1
        {
         X_\alpha \;\; \ar @<+2pt> `u[r] `[rr]^{g_{\alpha \beta}} [rr]
                                             \ar@{-->}[r]^{f_{\alpha}} 
         & \;\;X\;\;  \ar[r]^{h_{\beta}} 
         & \;\;Z  
         \\  
         S_\alpha \;\; \ar[r]^{\varphi_\alpha} 
         & \;\;S\;\;  \ar[r]^{\beta} 
         & \;\;u(Z)      
        }
$$
Consider the family of all arrows with domain $S$, $S\mr{\beta} u(Z)$, $Z \in \mathbb{T}$,   such 
that for all $\alpha$ there exists $X_\alpha\mr{g_{\alpha \beta}} Z$ in $\bb{T}$ over $\beta\circ\varphi_\alpha$. 
Take a initial family 
$X\mr{h_\beta} Z$ over the family $\beta$. For each $\alpha$, the family 
$g_{\alpha \beta}$ indexed by $\beta$ implies the existence of an arrow 
$X_\alpha\mr{f_\alpha} X$ 
over $\varphi_\alpha$ such that  $h_\beta\circ 
f_\alpha=g_{\alpha \beta}$. It remains to prove that 
$f_\alpha$ is a final family. Suppose we are given $Z_0 \in \bb{T}$, 
$S \mr{\beta_0} u(Z_0)$ in $\bb{S}$, and for all $\alpha$ an arrow 
$X_\alpha \mr{u_\alpha} Z_0$ over  $\beta_0 \circ \varphi_\alpha$. 
Let $h = h_{\beta_0}$. Since $u$ is faithful we have 
$u_\alpha = g_{\alpha \beta_0}$ so that $h \circ f_\alpha = u_\alpha$. Furthermore, again by the faithfulness of $u$, such an $h$ is unique since it must sit over $\beta_0$.    
\end{proof}

From the proof of this theorem we derive a further characterization of topological functors:
\begin{proposition} \label{chartop2}
A functor $u:\mathbb{T}\to\mathbb{S}$ is topological if and only if 

i) It is faithful and creates and preserves strict monomorphic
families.

ii) It has a full and faithful right adjoint  $u \dashv (-)_\top$,
$u(-)_\top = id$.
\end{proposition}
\begin{proof}
If $u$ is a topological functor, i) follows from proposition
\ref{topologicalisfaithful} and corollary \ref{fs=seintop}, and ii)
from corollary \ref{topologicalleftright}.
For the other direction, we shall show that $u$ creates
final families with the same proof that of proposition
\ref{topselfdual}. We show that the family in $\bb{S}$ over which an
initial family is created is actually a strict monomorphic
family. Consider the double diagram: 
$$
\xymatrix@1
        {
         X_\alpha \;\; \ar @<+2pt> `u[r] `[rr]^{g_{\alpha}} [rr]
                                             \ar@{-->}[r]^{f_{\alpha}} 
         & \;\;X\;\;  \ar[r]^{h_{id}} 
         & \;\;S_\top  
         \\  
         S_\alpha \;\; \ar[r]^{\varphi_\alpha} 
         & \;\;S\;\;  \ar[r]^{id} 
         & \;\;S      
        }
$$
where the arrows $g_\alpha$ correspond by adjointness to
$\varphi_\alpha$. This shows that $S \mr{id} S$ is in the family $S
\mr{\beta} u(Z)$ and thus
this family is strict monomorphic. By the dual of theorem
\ref{characterization2}, the strict monomorphic
family created over it is an initial family. The same proof
continuous then.  Notice that we are assuming $u$ to be faithful. 
\end{proof}

By proposition \ref{topselfdual} we have that topological functors are the same thing that functors satisfying the dual of definition \ref{topological}, that is, functors that create final families (or equivalently, cocartesian  families that compose). As a consequence, all the dual statements from \ref{topological} to \ref{chartop2} hold for topological functors. We explicitate the dual of this last characterization:

\begin{proposition} \label{chartop3}
A functor $u:\mathbb{T}\to\mathbb{S}$ is topological if and only if:

i) It is faithful and preserves and creates strict epimorphic
  families.

ii) It has a full and faithful left adjoint $(-)_\bot \dashv
\, u$, $u(-)_\bot = id$.
\end{proposition}

\begin{remark}
Notice that  condition i) in the previous propositions mean that a
topological functor is a $\cc{E}$-functor and a $\cc{M}$-functor.
\end{remark}

It is immediate to check that from the creation of strict epimorphic
(respectively,  strict monomorphic) families it follows the creation
of all colimits (respectively, limits) that may exists in $\bb{S}$
(see remark \ref{preservescolandlim}). We finish this section writing
down a list of properties that topological functors have.

\begin{theorem} \label{resumentopological}
If $u:\mathbb{T}\to\mathbb{S}$ is a topological functor (definition \ref{topological}), then:

1) It is faithful.

2) It creates final families and initial families.

3) It is a fibration and a cofibration, the fibers are complete (hence also cocomplete) lattices, and for any arrow $S \mr{\varphi} T$ in $\bb{S}$, there is a pair of adjoint functors  $\varphi\,! \dashv  \varphi^*$  between $\bb{T}_S$ and $\bb{T}_T$, 
$(-)^*$ is the action of the  fibration, and $!\,(-)$ is the action of the cofibration.

4) There are adjunctions $(-)_\top \vdash u \vdash (-)_\bot$, with $u(-)_\top = u(-)_\bot = id$.

5) In $\bb{T}$ strict epimorphic families are the same than final surjective families, and strict monomorphic families are the same than initial injective families.

6) It creates and preserve any limit and colimit that may exists in $\bb{S}$.

7) If in $\bb{S}$ all epis are strict (for example, when it is a topos), then in $\bb{T}$ epimorphic families are the same than surjective families.

\hspace{10pt} If in $\bb{S}$ all monos are strict (for example, when it is a topos), then in $\bb{T}$ monomorphic families are the same than injective families.

8) When $\bb{S} = \bb{S}et$ is the category of sets, then $u$ is representable by the object $1_\bot$
\end{theorem}

\vspace{2ex}

The forgetful functor from the category of topological spaces and of
the several quasitopoi associated with this category are examples of
topological functors (see \cite{D} and the references
therein). However, forgetful functors from concrete quasitopoi in a
more general sense are
not always topological functors. They lead to a notion of
quasitopological functor, which we shall develop elsewhere.

\end{document}